\documentclass[12pt,twoside]{article}
\usepackage{ifthen,latexsym,amssymb,amsmath}
\newtheorem{theorem}{Theorem}
\newtheorem{lemma}{Lemma}

{\begin{enumerate}}%
{\end{enumerate}}
\newcommand{\bull}{\mbox{\hspace{5 true pt}%
  \rule[-0.85 true pt]{3.9 true pt}{8.1 true pt}}}
\newcommand{\qed}{\mbox{\hspace{5 true pt}%
  \rule[-0.85 true pt]{2.0 true pt}{8.1 true pt}}}
\newenvironment{proof}{\par\smallbreak\noindent{\bf Proof.~}}%
{\unskip\nobreak\hfill \bull \par\medbreak}
\newenvironment{proofof}[1]{\par\smallbreak\noindent{\bf Proof~of~#1.~}}%
{\unskip\nobreak\hfill \bull \par\medbreak}
\newenvironment{subproof}{\par\noindent{\bf Proof of Claim.}}%
{\qed\par\smallbreak}
\newcounter{claim}
\renewcommand{\theclaim}{\arabic{claim}}
\newenvironment{claim}{\refstepcounter{claim}%
\par\medskip\par\noindent{\bf Claim~\theclaim.}\rm}%
{\par\medskip\par}
\newcommand{\setdef}[2]{\left\{ #1 \,:\, #2 \right\}}
\newcommand{\of}[1]{\left( #1 \right)}
\newcommand{\compl}[1]{\overline{#1}}

\newcommand{\rk}{\mathop{\mathit{rk}}\nolimits}
\newcommand{\envi}{\mathop{\mathit{Env}}\nolimits}
\newcommand{\env}[2]{\envi_{#1}(#2)}
\newcommand{\deco}{\mathop{\mathit{Dec}}\nolimits}
\newcommand{\dec}[2]{\deco_{#1} #2}
\newcommand{\case}[2]{\medskip\noindent\textbf{Case #1.} #2\medskip\par}
\newcommand{\subcase}[2]{\medskip\noindent\textbf{Subcase #1.} #2\medskip\par}
\newcommand{\tower}{\mathrm{Tower}}

\newcommand{\Alt}{\mathop{\sl alt}\nolimits}
\newcommand{\alt}[1]{\Alt(#1)}
\renewcommand{\mod}[1]{\,(\bmod\, #1)}            %(mod #1)
\newcommand{\ga}{\alpha}
\newcommand{\gb}{\beta}
\newcommand{\gs}{\sigma}

\newcommand{\EHRV}[1]{\mathit{Ehrv}(#1)}

\newif\ifnotesw\noteswtrue
\newcommand{\comment}[1]{\ifnotesw $\blacktriangleright$\ {\sf #1}\ 
  $\blacktriangleleft$ \fi}

\noteswfalse   % turn off marginal notes for now

\newcommand{\beq}[1]{\begin{equation}\label{eq:#1}}
\newcommand{\eeq}{\end{equation}}
\newcommand{\req}[1]{{\rm (\ref{eq:#1})}}

\newcommand{\args}[1]{{\rm[$#1$]}}
\newcommand{\Args}[1]{{\bf[$#1$]\,}}

\newcommand{\Ehr}{\mathrm{Ehr}}
\newcommand{\alternating}{alternating}

\newcommand{\KPSV}{kim+pikhurko+spencer+verbitsky:05}
\newcommand{\PSV}{pikhurko+spencer+verbitsky:05}
\newcommand{\PVV}{pikhurko+veith+verbitsky:04}

\title{
Decomposable graphs and\\
definitions with no quantifier alternation}

\author{Oleg Pikhurko\footnote{Partially supported by the National Science
Foundation,  Grant DMS-0457512.}\\
Department of Mathematical Sciences\\
Carnegie Mellon University\\ Pittsburgh, PA 15213, USA\\
Web: {\tt http://www.math.cmu.edu/\symbol{126}pikhurko}\\[2mm]
\and 
Joel Spencer\\
 Courant Institute\\ 
 New York University\\
  New York, NY 10012, USA
 \and
Oleg Verbitsky\footnote{Supported by an Alexander von Humboldt fellowship.}\\
Institut f\"ur Informatik\\
Humboldt Universit\"at Berlin\\
D-10099 Berlin, Germany}

\begin{document}

\maketitle

\begin{abstract}
Let $D(G)$ be the minimum quantifier depth of a first order sentence $\Phi$
that defines a graph $G$ up to isomorphism. Let $D_0(G)$ be the version 
of $D(G)$ where we do not allow quantifier alternations in $\Phi$. 
Define $q_0(n)$ to be the minimum of $D_0(G)$ over all graphs $G$ of order
$n$.

We prove that for all $n$ we have
 $$
 \log^*n-\log^*\log^*n-1\le q_0(n)\le \log^*n+22,
 $$
 where $\log^*n$ is equal to
the minimum number of iterations of the binary logarithm needed to bring $n$
to $1$ or below. 
The upper bound is obtained by constructing special graphs with
modular decomposition of very small depth.%
 \comment{I do not like calling is series-parallel decomposition, since
for me the name is strongly associated with the operations of
subdividing and duplicating edges from the construction of graphs without
$K_4$-minor.} 
 \end{abstract}


\comment{Some general remarks:

I decided to expand the introduction (using our paper in RSA) so that it
is understandable to mathematicians without any special background in
logic. On the other hand, I kept the description of the Ehrenfeucht game short
but gave references so that people who get interested can look up more
details if necessary.}

\section{Introduction}

We are interested in defining a given graph $G$ in first order logic, being
as succinct as possible. In order to state this problem formally, we have to
specify what we mean by the terms \emph{defining}, \emph{succinct}, etc. 

The vocabulary consists of the following symbols:
 \begin{itemize}
 \item \emph{variables} ($x$, $y$, $y_1$, etc);
 \item the \emph{relations} $=$ (equality) and $\sim$ (graph
adjacency);
  \item the \emph{quantifiers} $\forall$ (universality) and $\exists$
(existence);
 \item the usual Boolean \emph{connectives} ($\vee$,
$\wedge$, and $\neg$);
 \item  \emph{parentheses} (to indicate or change the precedence of operations). 
 \end{itemize}
 These can be combined into \emph{first order formulas} accordingly
to the standard rules. The term  \emph{first order} means that  the variables
represent vertices so the quantifiers 
apply to vertices only.%
 \comment{i.e.\ we cannot express properties like
\textsl{``There is a set $X$ having a given property''}. (In fact, the
language lacks any symbols to represent sets or functions.)\par}  
 In this paper, a \emph{sentence} is
a first order formula without free variables. On the
intuitive level it is perfectly clear what we mean when we say that a
sentence $\Phi$ \emph{is true} on a graph $G$. This is denoted by $G\models
\Phi$; we write $G\not\models \Phi$ for its negation ($\Phi$ \emph{is false} on
$G$).   We do not formalize these notions. A more detailed discussion
can be found in e.g.~\cite[Section~1]{spencer:slrg}.

Of course, if $G\models \Phi$ and $H\cong G$ (i.e.\ $H$ is isomorphic to
$G$), then $H\models \Phi$.  On the other hand, for any graph $G$ it is
possible to find a sentence $\Phi$ which \emph{defines} $G$, that
is, $G\models \Phi$ while $H\not\models \Phi$ for any $H\not\cong
G$. Indeed, let $V(G)=\{v_1,\dots,v_n\}$ be the vertex set of $G$ and $E(G)$
be its edge set. The required sentence
could read:
  \renewcommand{\arraystretch}{1.4}
 \beq{gen}
 \begin{array}{ccl}
 \Phi &=& \exists x_1\dots \exists x_n\ \left(\,\mathrm{Distinct}(x_1,\dots,x_n)\wedge
 \mathrm{Adj}(x_1,\dots,x_n)\,\right)\\
  &\wedge& \forall x_1\dots\forall x_{n+1}\ \neg\,
 \mathrm{Distinct}(x_1,\dots,x_{n+1}),\end{array}
 \eeq
 where, for the notational convenience, we use the following shorthands
 \begin{eqnarray*}
   \mathrm{Distinct}(x_1,\dots,x_k)&=& \bigwedge_{1\le i<j\le k} \neg\, (x_i=x_j)\\
 \mathrm{Adj}(x_1,\dots,x_n)&=&\bigwedge_{\{v_i,v_j\}\in E(G)} x_i\sim
 x_j\ \wedge\bigwedge_{\{v_i,v_j\}\not\in E(G)} \neg\, (x_i\sim
 x_j).
 \end{eqnarray*}
 In other words, we first specify that there are $n$ distinct vertices,
list the adjacencies and non-adjacencies between them, and then state that the total number of
vertices is at most~$n$.

A defining sentence $\Phi$ is not unique, so we are interested in finding
one which is
as succinct as possible. All natural succinctness measures of $\Phi$
are of interest: 
 \begin{itemize}
 \item the {\em length\/} $L(\Phi)$ which is the total number of
symbols in $\Phi$;
 \item the {\em quantifier depth\/} $D(\Phi)$ which is the maximum number 
of nested quantifiers in $\Phi$;
 \item the {\em width\/} $W(\Phi)$ which is the 
number of variables used in $\Phi$ (different occurrences of the same variable 
are not counted).
 \end{itemize}
 For example, for the sentence in~\req{gen} we have
$L(\Phi)=\Theta(n^2)$ and $D(\Phi)=W(\Phi)=n+1$. 
All three characteristics inherently arise in the analysis of
the computational problem of checking if a $\Phi$ is true on a given graph,
see e.g.\ Gr\"adel~\cite{gradel:05}.
They give us a small hierarchy of descriptive complexity measures
for graphs: $L(G)$ (resp.\ $D(G)$, $W(G)$) is the minimum of $L(\Phi)$ 
(resp.\ $D(\Phi)$, $W(\Phi)$) over all sentences $\Phi$ defining $G$.
These graph invariants will be referred to as
the \emph{logical length, depth}, and \emph{width} of $G$.
We have
 $$W(G)\le D(G)\le L(G).$$
The former number is of relevance for graph isomorphism testing, 
see Cai, F\"urer, and Immerman~\cite{cai+furer+immerman:92}.
The parameters 
$W(G)$ and $D(G)$ admit a purely combinatorial characterization in terms
of the Ehrenfeucht game, see \cite{cai+furer+immerman:92,spencer:slrg}.

Here, we address the logical depth of graphs which was recently studied
in Bohman et al~\cite{bohman+frieze+luczak+pikhurko+smyth+spencer+verbitsky:04,\KPSV,\PSV,\PVV,pikhurko+verbitsky:05,spencer+stjohn:05,verbitsky:05}. 
We focus on the following general question: How do restrictions on logic
affect the descriptive complexity of a graph?
Call a sentence $\Phi$ {\em $a$-\alternating\/} if
it contains negations only in front of relation symbols and every
sequence of nested quantifiers in $\Phi$ has at most $a$ quantifier
alternations, that is, the
occurrences of $\forall\exists$ and $\exists\forall$.
Let $D_a(G)$ denote the variant of $D(G)$ for $a$-\alternating\ defining
sentences. Clearly, for any integer $a\ge 0$ we have 
 $$
 D(G)\le D_{a+1}(G)\le D_a(G).
 $$

For example, the sentence in~\req{gen} has no alternations. Thus it shows
that for any graph $G$ we have 
 \beq{genD0}
 D_0(G)\le v(G)+1,
 \eeq
 where $v(G)$ 
denotes the number of
vertices in $G$. This bound is in general best possible: for example,
$D_0(K_n)=D(K_n)=n+1$. In Kim et al~\cite{\KPSV} we 
proved that $D(G)=\log_2n-\Theta(\log_2n\log_2n)$ and
$D_0(G)\le(2+o(1))\log_2n$ for almost all graphs $G$ of order $n$. 

In the above results, the functions $D(G)$ and $D_0(G)$ 
are the same or differ by at
most a constant factor. However, they can be very far apart in general. In
\cite[Corollary~5.7]{\PSV} we demonstrated a \emph{superrecursive gap} between
$D(G)$ and $D_0(G)$: 
namely, we proved that for any total recursive function $f$ there is a graph
$G$ with $f(D_0(G))< D(G)$. This is not too surprizing, since
the logic of 0-\alternating\ sentences is very restrictive
and provably weaker than the unbounded first order logic. Whereas 
the problem of deciding if a first order sentence is satisfiable
by some graph is unsolvable, 
it becomes solvable if restricted to 0-\alternating\
sentences. The last result is due to Ramsey's logical work \cite{ramsey:30}
founding the combinatorial Ramsey theory (see
Ne\v set\v ril~\cite[pp.~1336--1337]{nesetril:HC}
for historical comments on the relations between Ramsey theory and
logic).

Given Ramsey's decidability result, it is reasonable to concentrate on the
first order definability with no quantifier alternation. As our main result
here (Theorem~\ref{th:main}), we determine the asymptotic behavior of the
\emph{succinctness function $q_0(n)$}, where for an 
integer $a\ge 0$ we define
 $$
 q_a(n)=\min\setdef{D_a(G)}{G\mbox{\ has\ order\ }n}.
 $$
 Let \emph{log-star} $\log^*n$ be equal
to the minimum number of iterations of the binary logarithm needed to bring
$n$ to $1$ or below.

\begin{theorem}\label{th:main} For all $n$ we have
 \beq{main}
  \log^*n-\log^*\log^*n-1\le q_0(n)\le\log^*n+22.
 \eeq
 \end{theorem}

The estimates~\req{main} are 
in sharp contrast to the result in~\cite[Corollary~9.1]{\PSV} 
which shows a superrecursive gap between
 $$
 q(n)=\min\setdef{D(G)}{G\mbox{\ has\ order\ }n}
 $$
 and $n$. Thus Theorem~\ref{th:main}, 
besides being an interesting result on its own, implies that
we cannot have $q_0(n)\le f(q(n))$ for some total recursive $f$ and all $n$. 
This 
implies, again, a superrecursive gap between the graph invariants $D(G)$ 
and $D_0(G)$.

In \cite[Theorem~7.1]{\PSV} a weaker bound $q_0(n)\le 2\log^* n+O(1)$ for an infinite
sequence of values of $n$ is proved by inductively constructing large
asymmetric trees and estimating $D_0(G)$ in terms of their (very small)
radius. Here, our construction produces a graph of large order that has very
short \emph{modular decomposition\/} (as defined in Brandst\"adt, Le, and
Spinrad~\cite[Section~1.5]{brandstadt+le+spinrad:99}), starting with small
complement-connected graphs.  It seems feasible that many other recursively
defined constructions of graphs (see  Borie, Parker, and
Tovey~\cite{borie+parker+tovey:04} and Brandst\"adt, Le, and
Spinrad~\cite[Section~11]{brandstadt+le+spinrad:99} 
for surveys) may lead to upper bounds on $q_0(n)$
compatible with~\req{main}.  However, the proof of the upper bound
in~\req{main} required from
us many delicate auxiliary lemmas, even though we chose a construction which
is, in our opinion, most suitable for our purposes. So, a general theorem
would probably be very messy and difficult to prove.

In~\cite[Theorem~9.3]{\PSV} we have shown that 
 $$
 \log^*n-\log^*\log^* n-2\le q(n)\le \log^* n+4
 $$ 
 for infinitely many $n$. 
Combined with Theorem~\ref{th:main} and the obvious inequalities
 $$
 q_0(n)\ge q_a(n)\ge q(n),\quad \mbox{for any integer $a\ge 1$,}
 $$
 this implies that for any fixed $a$ we have $q_a(n)=(1+o(1))\log^* n$ for
infinitely many~$n$. We do not even know if $q_1(n)=(1+o(1))\log^* n$ for all
large $n$. 

In fact, Theorem~\ref{th:main} holds also for digraphs,
where instead of the adjacency relation $\sim$ we use the relation $x\mapsto
y$ to denote that the ordered pair $(x,y)$ is an arc. For example, the digraph
version of the lower bound in~\req{main} reads as follows.

\begin{theorem}\label{thm:lower:digraph} For
any digraph $G$ on $n$ vertices we have
 \beq{lower:digraph}
 D_0(G)\ge\log^*n-\log^*\log^*n-1.
 \eeq
\end{theorem}

Let us see how these results are related.  Take any graph $G$ and a
$0$-\alternating\ sentence $\Phi$ defining it. Let the digraph $G'$ be obtained from $G$ by
replacing each edge $\{x,y\}\in E(G)$ by a pair of arcs $(x,y)$ and
$(y,x)$. Then the sentence 
 $$
 \Big(\forall x\ \neg\, (x\mapsto x)\Big) \wedge\ \Big(\forall x\forall y\
\big((x\mapsto y)\wedge  (y\mapsto x)\big)\vee \big(\neg(x\mapsto y)\wedge  \neg 
(y\mapsto x)\big)  \Big) \wedge \Phi'
 $$
 defines $G'$, where $\Phi'$ is obtained from $\Phi$ by replacing each
occurence of $x\sim
y$ by, for example, $x\mapsto y$. Thus 
 $$
 D_0(G')\le \max(2,D_0(G))=D_0(G).$$
 This shows
that it is enough to prove the upper bound in Theorem~\ref{th:main} 
and the lower bound of Theorem~\ref{thm:lower:digraph}.

\section{Definitions}\label{s:prel}

We denote $[m,n]=\{m,m+1,\dots,n\}$ and $[n]=[1,n]$.
We define the \emph{tower-function} by $\tower(0)=1$ and
$\tower(i)=2^{\tower(i-1)}$ for each subsequent~$i$. Note that
$\log^*(\tower(i))=i$. The notation $x\in^i X$ means
$x\in X$ for odd $i$ and $x\not\in X$ for even $i$. (The mnemonic rule to
remember which is which is $\in^1=\in$.) 
The abbreviation `iff' means `if and
only if.' We do not allow infinite sentences nor infinite graphs (nor the
degenerate graph
with the empty vertex set).

We use the following graph notation: $\compl G$ is the complement of $G$;
$G\sqcup H$ is the vertex-disjoint union of graphs $G$ and $H$;
$G\subset H$ means that $G$ is isomorphic to an induced subgraph of $H$
(we will say that $G$ is {\em embeddable\/} into $H$).
For graphs (resp.\
sets) $A$ and $B$ the relation $A\subset B$ does not
exclude the case of isomorphism $A\cong B$ (resp.\ equality $A=B$).

We call $G$ {\em complement-connected\/} if both $G$ and $\compl G$
are connected. An inclusion-maximal complement-connected induced subgraph
of $G$ will be called a {\em comp\-lement-connected component\/} of $G$
or, for brevity, {\em cocomponent\/} of $G$. Cocomponents have no common 
vertices and their vertex sets partition $V(G)$.

The {\em decomposition\/} of $G$, denoted by $\deco G$, is the set of
all connected components of $G$ (this is a set of graphs, not just isomorphism 
types). Furthermore, given $i\ge0$, we define the
{\em depth $i$ decomposition\/} $\dec iG$ of $G$ by
 $$
 \dec 0G=\deco G\quad \mbox{and}\quad
\dec{i+1}G=\bigcup_{F\in\dec iG}\deco\compl F.
 $$
 Note that $\dec iG$ consists of connected graphs, and
distinct vertices $x,y$ of an $F\in \dec iG$ are adjacent in $F$ if and only if
$\{x,y\}\in^{i+1} E(G)$. Moreover,
 \beq{Pi}
 P_i=\setdef{V(F)}{F\in\dec iG}
 \eeq
 is a partition of $V(G)$
and $P_{i+1}$ refines $P_i$.
The {\em depth $i$ environment\/} of a vertex $v\in V(G)$,
denoted by $\env i{v;G}$, is the graph $F$ in $\dec iG$ containing~$v$. If the
underlying graph $G$ is clear from the context, we will usually write $\env
iv$.

We define the {\em rank\/} of a graph $G$, denoted by $\rk G$,
inductively as follows:
 \begin{itemize}
 \item If $G$ is complement-connected, then $\rk G=0$.
 \item 
If $G$ is connected but not complement-connected, then
$\rk G=\rk\compl G$.
 \item 
If $G$ is disconnected, then $\rk G=1+\max\setdef{\rk F}{F\in\deco G}$.
 \end{itemize}
Note that for connected graphs $\rk G$ is equal to the smallest $k$ 
such that $P_{k+1}=P_k$
or, equivalently, such that $P_k$ consists of $V(F)$ for all cocomponents $F$
of $G$.

Let $G$ be a connected graph and let $k=\rk G$. 
We call $G$ {\em uniform\/} if $\dec{k-1}G$
contains no complement-connected graph, that is, every cocomponent
appears in $\dec{k}G$ and no earlier. 
We call $G$ {\em inclusion-free\/} if the following two conditions
are true for every $0\le i\le k$: 
 \begin{enumerate}
 \item For any $K\in\dec iG$, $\compl K$ contains no isomorphic 
connected components.
  \item Of any two elements $K,M\in\dec iG$ none is properly
embeddable into the other, that is, either $K\cong M$ or none is an induced 
subgraph of
the other.
 \end{enumerate}

Let us now describe the {\em Ehrenfeucht game\/} $\Ehr_k(G,H)$ which will be
our tool
for studying the logical depth of graphs. The board consists of 
two vertex-disjoint graphs $G$ and $H$. There are $k$ rounds. The graphs $G,H$
and the number $k$ are known to both players, \emph{Spoiler} and
\emph{Duplicator} (or \emph{he} and \emph{she}).  In each round Spoiler 
selects one vertex in either $G$ or $H$; then Duplicator
must choose a vertex in the other graph. Let $x_i\in V(G)$ and $y_i\in V(H)$
denote the vertices selected by the players in the $i$-th round,
irrespectively of who selected them.
Duplicator wins the game if the componentwise correspondence between the
ordered $k$-tuples
$x_1,\ldots,x_k$ and $y_1,\ldots,y_k$ is a partial isomorphism from $G$ to
$H$. Otherwise the winner is Spoiler. In the {\em 0-alternation game\/}
Spoiler must play all the game in the same graph he selects in the first round.

Assume that $G\not\cong H$. Let $D(G,H)$ (resp.\ $D_0(G,H)$) denote the
minimum of
$D(\Phi)$ over all (resp.\ 0-\alternating) sentences $\Phi$ that are
true on one of the graphs and false on the other.  The Ehrenfeucht
theorem~\cite{ehrenfeucht:61} (see also Fra{\"\i}ss{\'e}~\cite{fraisse:54})
relates $D(G,H)$ and the length of the Ehrenfeucht game on $G$ and $H$. We
will use the following version of the theorem: $D_0(G,H)$ is equal to the
minimum $k$ such that Spoiler has a winning strategy in the $k$-round
0-alternation Ehrenfeucht game on $G$ and $H$. We will also use the fact
(see~\cite[Proposition~3.6]{\PSV}) that 
 $$D_0(G)=\max\setdef{D_0(G,H)}{H\not\cong G}.$$
 We refer the Reader to~\cite[Section~2]{spencer:slrg} which contains a
detailed discussion of the Ehrenfeucht game.

\section{Proof of the Upper Bound in Theorem \protect\ref{th:main}}

\subsection{Preliminaries}

\begin{lemma}\label{lm:cocosub}
Every complement-connected graph $G$ of order at least 5 has a vertex $v$ such 
that $G-v$ is still complement-connected.
\end{lemma}

\begin{proof} Suppose that the claim is false.
Take an arbitrary $v\in V$, where $V=V(G)$. This
vertex does not work so assume that, 
for example, $G-v$ is disconnected. Choose a
proper partition $V\setminus \{x\}=A_1\cup A_2$ such that no edge of $G$ connects $A_1$ to
$A_2$. Assume that $|A_1|\ge |A_2|$. Since $G$ is connected, the graph
$G_i=G[A_i\cup\{v\}]$ is connected, $i=1,2$. This implies that
$U_i\not=\emptyset$ for $i=1,2$, where 
 $$
 U_i=\{u\in A_i\mid G_i-u\mbox{ is connected}\}.
 $$
 Let $u\in U_1$. The graph $G-u$ is connected because any vertex of
$A_1\setminus\{u\}$  can be
connected (in $G_1-u$) to $v$  and then connected (in $G_2$) to any vertex of
$A_2$. 
Since $\compl G$ contains all edges between $A_1$ and $A_2$ (and $|A_1|\ge
2$), the graph $\compl
G-u-v$ is connected. Thus the only way that $u$ can fail to satisfy the
conclusion of the lemma is that $v$ is adjacent (in $G$) to every other vertex
except $u$ (the vertex $v$ cannot be adjacent to $u$ too because $G$ is complement-connected). 
The latter condition determines $u$ uniquely and therefore $U_1=\{u\}$. 
If $|A_2|\ge2$, then the same argument shows that $U_2$ should consist of the 
unique neighbor $u$ of $v$,
which is impossible. Thus, $|A_2|=1$ and hence $|A_1|\ge 3$. 
Let $w\in A_1$ be some neighbor of $u$
and let $z\in A_1\setminus\{u,w\}$. Then $G_1-z$ is still connected: $u$ is
connected to $v$ via $w$ while any other vertex is directly adjacent to
$v$. Hence, $z\in U_1$. This contradiction finishes the proof.
\end{proof}

\comment{
The following easy claim will be useful.

\begin{lemma}\label{lm:fact1} Let $H_0$ be a subgraph of some
member of $\dec l{H}$. Then 
for any $y,y'\in V(H_0)$ and an integer $i\ge0$, if $\env i{y;H_0}=\env i{y';H_0}$
then $\env {i+l}{y;H}=\env{i+l}{y';H}$.\hfill\bull
\end{lemma}
 }

Now we come to two strategic lemmas. The arguments of each lemma are listed
in square brackets. This is convenient when we refer back to these
results and, hopefully, makes the dependences between the lemmas 
easier to verify.

\begin{lemma}\Args{x,x',y,y',G,H,l}\label{lm:plung}
Consider the Ehrenfeucht game on graphs $G$ and $H$.
Let $x,x'\in V(G)$, $y,y'\in V(H)$ and assume that
the pairs $x,y$ and $x',y'$ were selected by the players in
the same rounds. Furthermore, assume that all the following properties hold.
 \begin{enumerate}
 \item $\env lx\ne\env l{x'}$.
 \item $\env ly=\env l{y'}$.
 \item $V(\env{l+1}{y})\not=V(\env{l}{y})$.
 \end{enumerate}
 Then Spoiler can win in at most $l+1$ rounds, playing all the time in $H$.
\end{lemma}

\begin{proof}
We proceed by induction on $l$. The induction step takes care of the base case
$l=0$ too. Observe that, for every $0\le i\le l$, we have 
$V(\env{i+1}{y})\not=V(\env{i}{y})$ so we do not have to worry about Assumption~3 when
using induction.

Let $m\in[0,l]$ be the minimum number such that $x'\notin\env mx$. If $m<l$,
Spoiler wins in $m+1\le l$ moves by induction.  So suppose that $m=l$.
Assume that $y$ and $y'$ are not adjacent in $\env ly$ for otherwise
Duplicator has already lost.  By Assumption~3 
the graph $\env{l}{y}$ is connected but not
complement-connected, so its diameter is at most $2$.  Spoiler selects any
$y''$ adjacent to both $y$ and $y'$ in $\env ly$. If Duplicator does not lose
in this round, it means that her reply $x''$ lies outside $\env{l-1}x$ (and
that $l\ge 1$). We have $\env{l-1}x\not=\env{l-1}{x''}$ and
$\env{l-1}y=\env{l-1}{y''}$. By the induction hypothesis applied to
\args{x,x'',y,y'',G,H,l-1}, Spoiler can win in at most $l$ extra
moves.\end{proof}

\begin{lemma}\Args{x_1,y_1,G,H,l}\label{lm:x1y1} Suppose that $x_1\in V(G)$
and $y_1\in V(H)$ were selected in some round of
the Ehrenfeucht game on $(G,H)$ so that there is an $l\ge 0$
satisfying the following \emph{Assumptions~1--3}.
 \begin{enumerate}
 \item $G_1=\env {l}{x_1}$ is not isomorphic to $H_1=\env {l}{y_1}$.
 \item $H_1$ is a uniform inclusion-free graph such that every
cocomponent of $H_1$ has at most $c$ vertices.
 \item For any $i\ge 0$, no member $A\in \dec i{H_1}$ is embeddable as a
proper subgraph  
into some $B\in \dec i{G_1}$.
 \end{enumerate}
 Then Spoiler can win the game in at most $k+c-1$ extra moves, playing
all the time inside $H$, where $k=\rk H_1+l$.\end{lemma}

\begin{proof} Suppose that it
is Spoiler's turn to move and, in addition to $x_1$ and $y_1$, 
we have the following configuration. Spoiler
has already selected vertices $y_2,\dots,y_s\in V(H_1)$, Duplicator has 
selected $x_2,\dots,x_s\in V(G_1)$, and all of the
following \emph{Properties~1--4} hold, where, for $j\in[s]$, we let
$H_j=\env{j+l-1}{y_j;H}$ and $G_j=\env{j+l-1}{x_j;G}$. 
 \begin{enumerate}
 \item For $i\in[2,s]$
we have $y_{i}\in V(H_{i-1})$.
 \item For $i\in[2,s]$ we have $x_{i}\in
V(G_{i-1})$.
 \item For every $i\in[s]$ we have
$H_i\not\cong G_i$.
 \item For every $i\in[2,s]$ the
vertices $y_i$ and $y_{i-1}$ belong to different components of
$\compl{H_{i-1}}$. (Note that $y_i\in V(H_{i-1})$ by Property~1.) 
 \end{enumerate}
 
Let us make a few remarks. Property~1 implies that
 $$
 V(H_1)\supset\ldots \supset V(H_s),
 $$
 and for $1\le i\le j\le s$ we have $y_j\in V(H_i)$.
Likewise by Property~2,
 \beq{Gnested}
 V(G_1)\supset \ldots\supset V(G_s),
 \eeq
 and for $1\le i\le j\le s$ we have 
$x_j\in V(G_i)$.  Properties~1 and~4 imply that $y_{j}\not\in V(H_i)$ for
any $1\le j<i\le s$. We stated Properties~1--4 this way in order 
to reduce the number
of checks needed to verify them. 
Also, note that we do not require that the vertices $x_i$ 
satisfy the analog of Property~4.

The above properties determine all $H$-adjacencies between the vertices
$y_1,\dots,y_s$. Indeed, take any $1\le i<j\le s$. By Properties~1 and~4,
$y_i$ and $y_j$ belong to different components of $\compl{H_i}$ so we have
$\{y_i,y_j\}\in E(H_i)$. This means that
$\{y_i,y_j\}\in^{i+l} E(H)$. In other words, 
the vertices $y_i$ and $y_j$ are adjacent in $H$ if and only if $i+l$ is
odd.

If $s=1$, then Properties~1, 2, and~4 are vacuously true, while Property~3 is
precisely Assumption~1 of the lemma.

We are going to show that Spoiler can either force the same situation after
the next round (of course, with $s$ increased by one) or win by making some
extra moves.

\case1{Suppose that $s< k-l$.}

As $H_s\not\cong G_s$, Assumption~3 (for $i=s-1$, $A=H_s$, and $B=G_s$) 
implies that $H_s\not\subset
G_s$. By Assumption~2, the connected graph 
$H_s\in\dec{s-1}{H_1}$ is inclusion-free; in particular, its complement
does not contain two isomorphic components. Hence, 
there is a component
$H_{s+1}$ of $\compl{H_s}$ which is not isomorphic to any component 
of~$\compl{G_s}$.

Suppose first that $y_s\not\in V(H_{s+1})$. Spoiler chooses an arbitrary
$y_{s+1}\in V(H_{s+1})$. Properties~1 and 4 hold automatically.  Let $x_{s+1}$
be Duplicator's reply. Assume that $x_{s+1}$ has the same adjacencies to the
previously selected vertices as $y_{s+1}$ for otherwise Spoiler has already
won having made $s\le k-l-1$ moves. (Note that we do not count $y_1$ as
a move, here or later in the proof.) 
Suppose that $x_{s+1}\not\in V(G_s)$, for otherwise Properties~2
and~3 hold and we are done.

\begin{claim}\label{cl:1}
 We have $l\ge 1$ and $x_{s+1}$ does not belong to $\env{l-1}{x_1;H}$.
 \end{claim}

\begin{subproof} First we argue that $x_{s+1}\not\in V(G_1)$. Suppose that
this is not true. In view of~\req{Gnested}, take the largest $i\in[s-1]$ such that $x_{s+1}\in V(G_i)$. By the
definition of $i$, $x_{s+1}\not\in V(G_{i+1})$, the latter being the
component of $\compl{G_i}$ that contains $x_{i+1}$. Thus
$\{x_{s+1},x_{i+1}\}\in E(G_i)$. 
 On the other hand,
$x_{s+1}$ is not adjacent to $x_{i+1}$  in $G_i$ because $y_{s+1}$ is not
adjacent to $y_{i+1}$ in $H_i$, a contradiction.

Next, we have $\{y_1,y_{s+1}\}\in^{l+1} E(H)$, so $\{x_1,x_{s+1}\}\in^{l+1}
E(G)$.  Since $x_{s+1}\not \in V(G_1)=V(\env l{x_1})$, we have $l\ge 1$. For
any vertex $z\in V(\env {l-1}{x_1})\setminus V(\env l{x_1})$ we have
$\{x_1,z\}\in^l E(G)$, so $x_{s+1}\not\in V(\env{l-1}{x_1})$, as
required.\end{subproof}

At this point it is possible to argue that, if $s\ge 2$, then Duplicator has
already lost. However, we still have to deal with the case $s=1$ (when we have
just $x_1$ and $x_2$). Since ruling out the case $s\ge 2$ would not make the
proof shorter, we do not do this.

We have $V(\env{l+1}{y_1})\not=V(\env{l}{y_1})$ because the latter set contains
$y_{s+1}$ while the former does not (or because $\rk H_1\ge l+s-1\ge 1$). 
Hence, Lemma~\ref{lm:plung} applies to
\args{x_1,x_{s+1},y_1,y_{s+1},G,H,l-1}, and Spoiler can win the game in at
most $l$ extra moves, having made at most $s + l \le k-1$ moves in total.

It remains to describe Spoiler's strategy if $y_s\in V(H_{s+1})$, when
Spoiler cannot just choose some $y_{s+1}\in V(H_{s+1})$ as this would
violate Property~4. Here, Spoiler first selects
some $y_{s+1}\in V(H_s)\setminus V(H_{s+1})$. (This set is non-empty since
$s<\rk H_1$.) Let Duplicator reply with
$x_{s+1}$. If $x_{s+1}\not\in V(G_s)$, then by the argument of 
Claim~\ref{cl:1} we have that $l\ge1$ and $\env
{l-1}{x_1}\not=\env{l-1}{x_{s+1}}$. Thus 
Spoiler can win in at most $l$ further
moves by Lemma~\ref{lm:plung},
having made at most $s + l \le k-1$ moves in
total. Hence, let us assume that $x_{s+1}\in V(G_s)$. In this case, let
us swap the vertices $y_s$ and $y_{s+1}$ as well as $x_s$ and
$x_{s+1}$. It is clear that the new sequences $y_1,\dots,y_{s+1}$
and $x_1,\dots,x_{s+1}$ satisfy Properties~1--4.
This completes the description of
the case $s<k-l$.

\case2{Suppose that $s=k-l$ (or that $s=1$ and $k=l$).}

This means that $H_s$ is a cocomponent of $H_1$ (and thus has at most $c$
vertices).  Spoiler selects all vertices in $V(H_s)\setminus\{y_s\}$. We claim
that Duplicator has lost by now. Indeed, if
Duplicator replies all the time inside $G_s$, then she has lost because
$G_s\not\supset H_s$ by
Assumption~3 and Property~3. Otherwise, her response to the whole set $V(H_s)$
cannot be complement-connected because it contains both a vertex outside of
$G_s$ and the vertex $x_s\in V(G_s)$. Thus Spoiler wins, 
having made at most $s-1+c-1\le k+c-1$ further moves.\end{proof}

\subsection{Finishing the Proof}

\begin{lemma}[Main Lemma]\label{lm:main2} Let $G$ be a connected 
uniform inclu\-sion-free graph.  Let $c\ge 5$ and 
suppose that every cocomponent of $G$ has
at most $c$ vertices. Then $D_0(G)\le \rk G+c+1$.
 \end{lemma}

%\vspace{-3mm}

\begin{proof} 
Let $k=\rk G$. Since the case of $k=0$ is trivial (namely we have $D_0(G)\le
v(G)+1\le c+1$ by~\req{genD0}),
we  assume that $k\ge1$.

Fix a graph $H\not\cong G$.  We will design a strategy allowing Spoiler to
win the 0-alternation Ehrenfeucht game on $(G,H)$ in at most the required
number of moves. There are a few cases to consider.

\case 1{$H$ has a cocomponent $C$ non-embeddable into any cocomponent of~$G$.}

If $C$ has no more than $c$ vertices,
Spoiler selects all vertices of $C$. Otherwise he selects $c+1$ vertices spanning 
a complement-connected subgraph in $C$ which is possible by
Lemma~\ref{lm:cocosub} (since $c\ge 5$).  If Duplicator's response $A$
is within a cocomponent of $G$, then $C\not\cong A$ by the assumption. 
Otherwise $A$ is not complement-connected and Duplicator loses anyway.

\case 2{There are an $l\in[0,k]$ and an $A\in\dec lG$
properly embeddable into some $B\in\dec lH$, and not Case 1.}

Let $H_0$ be a copy of $A$ in $B$. Fix an arbitrary vertex $y_0\in
V(B)\setminus V(H_0)$. Note that since we are not in Case 1, the
connected graph $B$ cannot be a cocomponent of $H$ by Property~2 in the
definition of an inclusion-free graph. Hence
 \beq{plung_y0}
 V(\env l{y_0;H})\not= V(\env{l+1}{y_0;H}).  
 \eeq
 
Let $Z=V(B)\setminus V(H_0)$. We will need the following routine claim, whose
proof uses~\req{plung_y0} and the connectedness of $H_0$.

\begin{claim}\label{cl:H-Z} 
 For any $m\ge 0$ and $y\in V(H_0)$ we have 
 $$\env m{y;H_0}=\env {m+l}{y;H-Z}.$$
\end{claim}
\begin{subproof} It is enough to prove the case $m=0$ only,
because the remaining cases would follow by a straightforward induction on
$m$. Since $H_0$ is connected, the claim for $m=0$ amounts to proving that 
$H_0=\env {l}{y;H-Z}$. The latter identity is precisely the case $s=l$ of 
 \beq{ind}
 \env{s}{y;H-Z}=\env{s}{y;H}-Z,\quad \mbox{for any $s\in[0,l]$}.
 \eeq

We prove~\req{ind} by induction on $s$. 
the case $s=0$ being routine to check. Let
$s\in [l]$. By~\req{plung_y0} the complement of $\env{s-1}{y;H}$ has at least
two components, one of which, namely $\env{s}{y;H}$, contains $V(H_0)$. In
order to prove~\req{ind} by induction on $s$ we
have to show that  $F=\env{s}{y;H}-Z$ is still connected. If $s=l$, then this is
true because $F=H_0$.  If $s<l$, 
then $F$ is connected because it 
contain a spanning complete bipartite graph with one part being
$V_1=V(\env{s+1}{y;H})\setminus Z$. (This bipartite graph is not degenerate:
$V_1\supset V(H_0)\not=\emptyset$
while  $V_1\not=V(F)$ by~\req{plung_y0}.)\end{subproof}

Spoiler plays in $H$.  At the first move he selects $y_0$. Denote Duplicator's
response in $G$ by $x_0$ and set $G_0=\env l{x_0}$. There are two alternatives
to consider.

\subcase{2.1}{$G_0\not\cong H_0$.}

Suppose first that $l<k$. Since $G_0$ and $H_0$ are non-isomorphic copies
of elements of $\dec lG$ and $G$ is inclusion-free, Spoiler is able to
make his next choice $y_1$ in some $H_1\in\deco\compl{H_0}$ with no isomorphic
graph in 
$\deco\compl{G_0}$.  Denote Duplicator's response by~$x_1$.

If $x_1\not\in V(G_0)$, then
Lemma~\ref{lm:plung} applies to \args{x_0,x_1,y_0,y_1,G,H,l} in view
of~\req{plung_y0}.  Thus 
Spoiler can win by using at most $l+3\le k+2$
moves in total. So, assume that
$x_1\in V(G_0)$. 
Lemma~\ref{lm:x1y1} applies to \args{x_1,y_1,G,H-Z,l+1} in view of
Claim~\ref{cl:H-Z}. (For example, Assumption~3 is satisfied because $G$ is
uniform inclusion-free and 
$\dec lG$ contains both $G_0$ and an
isomorphic copy of $H_0$.) Thus Spoiler can win in at most $2+(k+c-1)$ moves in
total, as desired.

It remains to consider the case $l=k$. Spoiler selects all vertices of $H_0$.
There are at most $c$ of them because $H_0$ is isomorphic to a cocomponent of
$G$. If Duplicator's replies lie in $V(G_0)$, she has already lost in view of
$G_0\not\supset H_0$ (which holds since $G$ is inclusion-free). Otherwise,
Duplicator's reply to $V(H_0)$ contains both a vertex outside $G_0$ and the
vertex $x_0\in V(G_0)$, so it cannot be complement-connected, and she loses.
So, Spoiler wins having made at most $c$ moves in total.

\subcase{2.2}{$G_0\cong H_0$.}  

Though the graphs are isomorphic, the crucial fact is that $G_0$, unlike
$H_0$, contains a selected vertex. By the definition of an inclusion-free
graph, every automorphism of $G_0\cong H_0$ takes each cocomponent onto
itself.  Therefore all isomorphisms between $G_0$ and $H_0$ match
cocomponents of these graphs in the same way. Let $Y$ be the $H_0$-counterpart
of the cocomponent $X=\env {k-l}{x_0;G_0}$ with respect to this matching. In
the second round Spoiler selects an arbitrary $y_1$ in $Y$. Denote
Duplicator's answer by $x_1$.

Suppose first that $x_1\in X$.  
Spoiler selects all 
vertices  of $Y\setminus\{y_1\}$.  At least one of
Duplicator's replies lies outside $V(X)$ for otherwise she has already lost
having chosen some vertex in $X$ twice. But then Duplicator's reply to $Y$
cannot be complement-connected. In any case Spoiler wins,
having made at most $c+1$ moves in total.

If $x_1\in V(G_0)\setminus X$, then there is an $m\le k-l$ such
that $\env m{x_1;G_0}$ and $\env m{y_1;H_0}$ are
non-isomorphic. By Claim~\ref{cl:H-Z} Spoiler can apply the strategy of
Lemma~\ref{lm:x1y1} to \args{x_1,y_1,G,H-Z,l+m}, winning in at most
$2+(k+c-1)$ moves. 
If $x_1\not\in V(G_0)$, then Spoiler wins by Lemma~\ref{lm:plung} applied to
\args{x_0,x_1,y_0,y_1,G,H,l}, having made at most $2+l+1< k+c+1$ moves in
total.

\case3{$H$ has a component $H_0$ isomorphic to $G$, and not Cases~1--2.}

Spoiler plays in $H$. In the first round he selects a vertex $y_0$ outside
$H_0$ and further plays exactly as in Subcase 2.2 with $G_0=G$.

\case 4{Neither of Cases 1--3.}
 
Spoiler plays in $G_0=G$. His first move $x_0$ is arbitrary.  Denote
Duplicator's response in $H$ by $y_0$ and set $H_0=\env0{y_0}$.
Since we are not in Cases 1--3, 
$G_0\not\subset H_0$.  As $G_0$ is inclusion-free,
$\compl{G_0}$ has a connected component $G_1$ with no isomorphic component in
$\compl{H_0}$. 

If $x_0\not\in V(G_1)$, then Spoiler just selects any vertex $x_1\in V(G_1)$.
Let Duplicator respond with $y_1$. Assume that $y_1\in V(H_0)$, for otherwise
Duplicator has already lost: $\{y_0,y_1\}\not\in E(H)$ while $\{x_0,x_1\}\in
E(G)$.

If $x_0\in V(G_1)$, then Spoiler selects any vertex $x_1\in V(G_0)\setminus
V(G_1)$. (The latter set is non-empty since $k\ge 1$.) 
Let Duplicator respond with $y_1$. As before we can assume that $y_1\in
V(H_0)$. Now, let us swap $x_0$ and $x_1$ as well as $y_0$ and $y_1$.
 
What we have achieved in both cases is that $G_1\not\cong H_1$, where
$H_1=\env1{y_1;H}$. Also, $G_1$ is a uniform
inclusion-free graph of rank $k-1$.
Lemma~\ref{lm:x1y1} applies to \args{y_1,x_1,H,G,1}. (For example, Assumption~3
of the lemma holds because we are not in Cases~1--2.) This
shows that Spoiler can win the $0$-alternation game in at most $2+ (k+c-1)=
k+c+1$ moves. This completes the proof of Lemma~\ref{lm:main2}. 
\end{proof}

\noindent\textbf{Proof of the upper bound in Theorem~\ref{th:main}.}
Fix an integer $c$ so that there are $4c+4$ pairwise non-embeddable into each
other complement-connected graphs 
 $$
 H_{i,j},\quad c\le i\le 2c,\ 1\le j\le4,
 $$
 such that $H_{i,j}$ has order $i$. The existence of $c$ can be easily deduced by
choosing each $H_{i,j}$ uniformly at random from all graphs of order $i$,
independently from the other graphs.
Indeed, for any $c\le i\le f\le 2c$ with $(i,j)\not=(f,g)$ the probability
that $H_{ij}$ is embeddable into $H_{fg}$ is at most 
 $$
\frac{f!}{(f-i)!}\, 2^{-{i\choose 2}}$$
 while the probability of $H_{i,j}$ not being
complement-connected is at most 
 $$\frac12\, \sum_{h=1}^{i-1} {i\choose h} 2^{-h(i-h)+1},$$
 where the factor $\frac12$ acounts for the fact that each vertex partition
is counted twice.

 Hence, by looking at the expected number of `bad' events, we conclude that if
 \beq{c}
  16 \sum_{c\le i\le f\le 2c} \frac{f!}{(f-i)!}\, 2^{-{i\choose 2}} + 
 4\times \frac12 \sum_{i=c}^{2c} \sum_{h=1}^{i-1} {i\choose h} 2^{-h(i-h)+1}<1,
 \eeq
 then the required graphs exist. The exact-arithmetic calculation with
\emph{Mathematica} shows that  $c=10$
works in~\req{c}. (This value of $c$ can perhaps be improved with more work.)%
 \comment{See the latex source for 
the estimate of the expected number of bad events I used with Mathematica. 
% Table[{c,
%    16*Sum[Factorial[f]/(2^Binomial[i, 2]*Factorial[f - i]),
%          {f, c, 2c}, {i, c, f}] + 
%      4*(1/2)*Sum[
%          Binomial[v, h]*2^(-h(v - h) + 1), {v, c, 2c}, {h, 1, v - 1}] - 
%      1}, {c, 9, 10}]
 }

We define, inductively on $i$, a family $R_i$ of graphs, starting with
 $$
 R_0=\{H_{c,1},H_{c,2},H_{c,3},H_{c,4}\}.
 $$
 Assume that $R_{i-1}$ is already specified.
Given a non-empty subset $S\subset R_{i-1}$, we define the graph 
 $$G_{i,S}=\compl{\bigsqcup_{G\in S} G},$$
 or, in words, $G_{i,S}$ is the complement of the vertex-disjoint union of
the graphs in $S$. We let
 $$
 R_{i}=\setdef{G_{i,S}}{|S|=|R_{i-1}|/2},
 $$
 where we view $R_i$ as 
the set of isomorphism types of graphs. It is proved in 
Claim~\ref{cl:Gi} below 
that the graphs $G_{i,S}$ are pairwise non-isomorphic. (In particular, 
this implies 
by induction on $i$ that $|R_i|$ is even because
${2m\choose m}$ is even for any integer $m\ge 1$.) Let $r_i=|R_i|$.

Let us list some properties of these graphs.

\begin{claim}\label{cl:Gi}
\begin{enumerate}
 \item
 For any $S\subset R_{i-1}$ with $|S|\ge 2$, 
$G_{i,S}$ is a connected inclusion-free uniform graph of rank $i$.
 \item For any $S,T\subset R_{i-1}$ with $S\not\subset T$, the graph
$G_{i,S}$ is not embeddable into $G_{i,T}$.
 \item $r_i={r_{i-1}\choose r_{i-1}/2}$.
\end{enumerate}\end{claim}
\begin{subproof}
 We prove all claims
by induction on $i$, the case $i=1$ directly following from the
definition of $R_0$.  Let $i\ge 2$. 

First, we verify Property~1, assuming that Properties~1--3 hold for all smaller
values of $i$. 
Since $|S|\ge 2$, $G_{i,S}$ is connected.  The components of $\compl{G_{i,S}}$ 
belong to $R_{i-1}$, each being isomorphic to
$G_{i-1,S'}$ for some $S'\subset R_{i-2}$. From Property~3 and the initial 
value $r_0=4$, it is easy to deduce that 
$|S'|=r_{i-2}/2\ge 2$. By the inductive Property~1, all
components of $\compl{G_{i,S}}$
are uniform of rank $i-1$, so $G_{i,S}$ is uniform of rank $i$. 

Next, let us verify that $G_{i,S}$ is an inclusion-free graph.  For any
$j\in[i]$ all elements of $\dec j{G_{i,S}}$ belong to $R_{i-j}$; by induction,
each is inclusion-free. Let us show
that none of these graphs is properly embeddable into another.%
 \comment{We have to prove this: even if each component of $\compl G$ is
inclusion-free, it may happen in principle that of
two elements of $\dec jG$ one is embeddable into another.}
 Assume that
$j<i$ for otherwise the claim follow from the definition of $R_0$. Take
any two non-isomorphic $G_{i-j,S'}, G_{i-j,S''}\in R_{i-j}$. 
We have $S'\not\subset S''$ because $S'\not=S''$ and
$|S'|=|S''|=r_{i-j-1}/2$. By induction (Property~2), we conclude that
$G_{i-j,S'}\not\subset  G_{i-j,S''}$, giving the stated. 
Since $G_{i,S}$ is connected, it remains to observe that
$\compl{G_{i,S}}$ has no two isomorphic components, which follows from
Property~2 again. Thus $G_{i,S}$ is indeed inclusion-free. We have
completely finished the inductive step for Property~1. 

Let us turn to Property~2. 
All components of $\compl{G_{i,S}}$ and $\compl{G_{i,T}}$ belong to 
$R_{i-1}$. Take any $H \in S\setminus T$. The graph $H\in
R_{i-1}$ appears as a component in $\compl{G_{i,S}}$. By induction
(Property~2) and the definition of $R_{i-1}$,
$H$ cannot be embedded into any component of $\compl{G_{i,T}}$. Thus
$G_{i,S}\not\subset G_{i,T}$, as required. Property~3 follows
from Property~2 which implies that the graphs $G_{i,S}$, for $S\subset R_{i-1}$,
are pairwise non-isomorphic.\end{subproof}

All graphs in $R_i$ have the same order which we denote by $n_i$. We have
$n_0=c$ and, for $i\ge 1$,
 $$
 n_i=n_{i-1}r_{i-1}/2.$$

 We have $v(G_{i,S})=|S|\,n_{i-1}$. If we denote $m_i=r_i/2$, then we have $m_0=2$ and $m_1=3$. Thus
for $i\ge 1$ we have
 $$
 m_{i+1}= \frac12\, r_{i+1} = \frac12\, {r_i\choose r_i/2} = \frac12\,
{2m_i\choose m_i} \ge 2^{m_i}.
 $$
 We conclude that $m_i> \tower(i)$ for all $i\ge 0$ and thus 
 \beq{ni}
 n_i\ge m_{i-1}> \tower(i-1).
 \eeq

At this point we are able to prove the required upper bound on $q_0(n)$
for an infinite sequence of $n$, namely,
 \begin{equation}\label{eq:Ni}
 \ldots, n_{i-1}, 2n_{i-1}, 3n_{i-1}, \ldots, m_{i-1}n_{i-1}=n_i, 
2n_i,\ldots 
 \end{equation}
 Indeed, by Lemma~\ref{lm:main2} for every $2\le s\le m_i$ and an
$s$-set $S\subset R_{i-1}$, we have 
 $$
  q_0(sn_{i-1}) \le D_0(G_{i,S})\le  i+c+1.
  $$
 Also, we have $i\le\log^*n_{i-1}+1$ by~\req{ni}. Thus
 $$
 q_0(sn_{i-1})\le \log^*(n_{i-1})+c+2 \le \log^*(sn_{i-1})+12.
 $$
 
It now remains to fill in the gaps in \req{Ni}. We need some auxiliary
notions and claims first. 
We define the operation of a \emph{cocomponent replacement} as follows.
Suppose that $A$ is a cocomponent of a graph $G$ and $B$ is a
complement-connected 
graph. The result of the \emph{replacement of $A$ with $B$ in $G$} is the graph $G'$ with
$V(G')=(V(G)\setminus V(A))\cup V(B)$ such that $G'[V(B)]=B$, $G'-B=G-A$, and
every vertex in $B$ is adjacent to a vertex $v$ outside $B$ in $G'$ if and
only if
every vertex in $A$ is adjacent to $v$ in $G$. (Here, we assume that $V(G)\cap
V(B)=\emptyset$, and we use the fact that any two vertices inside a cocomponent have
the same adjacency pattern to the rest of the graph.)

\begin{claim}\label{cl:spusk}
Let $G$ be a uniform inclusion-free graph of rank $i$ 
with all cocomponents being isomorphic to one of $H_{c,l}$ with $1\le l\le4$.
Let $G'$ be obtained from $G$ by replacing each cocomponent $A\cong
H_{c,l}$ with some $H_{j,l}$, where $j\in[c,2c]$ may depend on $A$. 
Then $G'$ is a uniform inclusion-free graph of rank $i$.
\end{claim}

\begin{subproof}
The partitions $P_0,\dots,P_i$ defined in~\req{Pi} are completely determined
by the vertex sets of the 
cocomponents and the adjacencies between then. This
shows that $G'$ is uniform of rank $i$. Let us check that $G'$ is
inclusion-free. 

Let $0\le j\le i$, $K'\in \dec j{G'}$, and $C_1',C_2'$ be some distinct 
components
of the complement of $K'$. Suppose on the contrary that a bijection 
$f':V(C_1')\to V(C_2')$ establishes an isomorphism between $C_1'$ and
$C_2'$. The isomorphism $f'$ induces a correspondence $g'$ between the
cocomponents of $C_1'$ and $C_2'$.

The definition of component replacement allows us
to point the corresponding $K\in \dec jG$, $C_1,G_2\in \deco{\compl
K}$, and $g$. Since $C_1\not\cong C_2$, there is a cocomponent $X_1$ 
of $C_1$ such that the cocomponent $X_2=g(X_1)$ is not isomorphic to $X_1$. It
means that, if $X_1\cong H_{c,l_1}$ and $X_2\cong H_{c,l_2}$, then
$l_1\not=l_2$. But in $G'$ these are replaced by $X_1'\cong H_{j_1,l_1}$ and
$X_2'\cong H_{j_2,l_2}$, which are still non-isomorphic since
$l_1\not=l_2$. This contradicts the assumption that $f'$ is an
isomorphism. Thus $G'$ satisfies Property~1 
of the definition of an inclusion-free
graph. The other property in the definition can be checked similarly.
\end{subproof}

If $n\le 2c=20$, then the upper bound~\req{main} follows from the trivial
inequality $q_0(n)\le n+1$. So assume that $n> 2c=2n_0$.
Choose the integer $i$ satisfying
$2n_i\le n < 2n_{i+1}$. Since $n_{i+1}=n_im_i$, let $s\in
[2,2m_i-1]$ satisfy $sn_i\le n< (s+1)n_i$. Pick any $s$-set $S\subset R_{i}$
and let $G=G_{i+1,S}$. 
We have $v(G)=sn_i\le n$ and, by Claim~\ref{cl:Gi}, the graph 
$G$ is inclusion-free and
uniform of rank $i+1$.

Let $f:\dec{i+1}{G}\to [c,2c]$ be some function. 
We construct a new graph $G_f$ by replacing 
every cocomponent $A$ of $G$ by a copy of $H_{f(A),j}$, where $j$
is defined by $A\cong H_{c,j}$. If $f$ is the constant function assuming the
value $2c$,  then $v(G_f)=2v(G)>n$. Hence there is some choice of $f$ such
that $v(G_f)=n$. By Claims~\ref{cl:Gi} and~\ref{cl:spusk}, the graph
$G_f$ is a uniform
inclusion-free graph of rank $i+1$. By Lemma~\ref{lm:main2}, we have $D_0(G_f)\le
i+2c+2$. On the other hand, $n\ge n_{i}> \tower(i-1)$, that is, $\log^*
n\ge i$. It means that
 $$
 q_0(n)\le \log^*n +2c+ 2=\log^*n+22.
 $$
 This finishes the proof of the upper bound in~\req{main}.\bull

\section{Lower bound: Proof of Theorem \protect\ref{thm:lower:digraph}}

From now on we will be dealing with digraphs.

Given a first order formula $\Phi$ in which the negation sign occurs only in
front of atomic subformulas, let
the {\em alternation number\/} of ${\Phi}$, denoted by $\alt\Phi$,
be the maximum number of quantifier alternations, i.e.\
the occurrences of $\exists\forall$ and $\forall\exists$, in
a sequence of nested quantifiers of~$\Phi$. 
For a non-negative integer
$a$, we denote 
 $$
 \Lambda_a=\setdef{\Phi}{\alt\Phi\le a}.
 $$
 We also define $\Lambda_{1/2}$ to be the class of formulas $\Phi$ with $\alt\Phi\le1$ such that
any sequence of nested quantifiers of $\Phi$ starts with $\exists$ or has
no quantifier alternation. Note that
$\Lambda_0\subset\Lambda_{1/2}\subset\Lambda_1\subset\Lambda_2\subset\ldots$.

Now we somewhat extend our notation.
Let $F$ be some class of first order formulas.
If a digraph $G$ has a defining sentence in $F$,
let $D_F(G)$ (resp.\ $L_F(G)$) denote the minimum quantifier rank
(resp.\ length) of a such sentence; otherwise, we let $D_F(G)=L_F(G)=\infty$.
The {\em succinctness function\/} is defined as
$$
q_F(n)=\min\setdef{D_F(G)}{v(G)=n}.
$$
 Whenever the index $F$ is omitted, it is supposed that $F$ is the class
of all first order formulas.
We also simplify notation by $D_a(G)=D_{\Lambda_a}(G)$ and similarly with
$L_a(G)$ and $q_a(n)$. Clearly,
 $$
 q(n)\le\ldots\le q_2(n)\le q_1(n)\le q_{1/2}(n)\le q_0(n).
 $$

\begin{lemma}\label{thm:lvsda}  For every $a\in\{0,1/2,1,2,3,\dots\}$ and any
digraph $G$ we have
 $$
 L_a(G) < \tower(D_a(G)+\log^*D_a(G)+2).
 $$
\end{lemma}

\noindent
An analog of this lemma for $L(G)$ and $D(G)$ appears in
\cite[Theorem 10.1]{\PSV}. However, the proof of Lemma \ref{thm:lvsda}
we give below is not just an easy adaptation of the proof in \cite{\PSV}
because the restrictions
on the class of formulas do not allow to run the same argument directly.
Moreover, if $a=1/2$, there appears another 
obstacle --- the class of formulas $\Lambda_{1/2}$
is not closed with respect to negation.

Lemma \ref{thm:lvsda} is proved in the next section in a stronger
form since the argument is presentable more naturally
in a more general situation. Here, let us show how Lemma~\ref{thm:lvsda}
implies Theorem \ref{thm:lower:digraph}.

Given $n$, denote $k=q_0(n)$ and fix a digraph $G$ on $n$ vertices such that
$D_0(G)=k$. By Lemma \ref{thm:lvsda}, $G$ is definable by a 0-\alternating\ sentence
$\Phi$ of length less than $\tower(k+\log^*k+2)$. First,
we convert $\Phi$ to an equivalent \emph{prenex $\exists^*\forall^*$-sentence}
$\Psi$, i.e.\ of form~\req{Psi}. 
This can be easily done as follows. By renaming variables,
ensure that each variable is quantified exactly once. Let the existential
(resp.\ universal) quantifiers appear with variables
$x_1,\dots,x_l$ (resp.\ $y_1,\dots,y_m$) in this order as we scan $\Phi$ from
left to right. To obtain the required sentence $\Psi$ simply `pull'
all quantifiers at front:
 \beq{Psi}
 \Psi=\exists x_1\, \dots\,\exists x_l\, \forall y_1\,\dots\, \forall y_m\
 (\mbox{quantifier-free part})
 \eeq
 The obtained sentence $\Psi$ is equivalent to $\Phi$ since the latter does
not contain an $\exists$-quantifier in the range of a
$\forall$-quantifier. Also, this reduction does not increase the total
number of quantifiers. Therefore, as a rather rough estimate, 
we have $D(\Psi)\le L(\Phi)$. 

It is well known and easy to see that, if a sentence of the form~\req{Psi}
is true on some structure $H$, then it is true on some structure of order
at most $l\le D(\Psi)$. (Indeed, fix any satisfying assignment for
$x_1,\dots,x_l$ and take the substructure of $H$ induced by the corresponding
vertices.)  Since the defining sentence $\Psi$
is true only on $G$, we have 
 $$
 n\le D(\Psi)\le L(\Phi)< \tower(k+\log^*k+2).
 $$
 This implies that 
 \beq{n1}
 \log^* n\le k+\log^*k+1.
 \eeq
 Suppose on the contrary to Theorem~\ref{thm:lower:digraph} that $k\le
\log^*n-\log^*\log^*n-2$. Then $\log^*k\le \log^*\log^*n$ and~\req{n1} implies
that 
 $$
 \log^* n \le (\log^*n-\log^*\log^*n-2) + \log^*\log^*n+1,
 $$
 which is a contradiction, proving Theorem~\ref{thm:lower:digraph}.

Note that identically the same argument works for $a=1/2$ as well, giving that
 \beq{q12}
 q_{1/2}(n)\ge\log^*n-\log^*\log^*n-1\quad\mbox{ for all $n$}.
 \eeq

\section{Length vs.\ depth for restricted classes of defining sentences}%
\label{s:lvsd}

Writing $A(x_1,\ldots,x_s)$, we mean that $x_1,\ldots,x_s$ are all
free variables of $A$. We allow $s=0$ which means that $A$ is a sentence.
A formula $A(x_1,\ldots,x_s)$ of quantifier rank $k-s$ is {\em normal\/} if
\begin{itemize}
\item
all negations occurring in $A$ stay only in front of atomic subformulas,
\item
$A$ has occurrences of variables $x_1,\ldots,x_k$ only,
\item
every sequence of nested quantifiers of $A$ has length $k-s$ and
quantifies the variables $x_{s+1},\ldots,x_k$ exactly in this order.
\end{itemize}
A simple inductive syntactic argument shows that any
$A(x_1,\ldots,x_s)$ has an equivalent normal formula $A'(x_1,\ldots,x_s)$
of the same quantifier rank. Such a formula $A'$ will be called
a {\em normal form\/} of~$A$.

Recall that by $F$ we denote a class of first order formulas.
Given $F$, the class of sentences (i.e.\ closed formulas)
in $F$ of quantifier rank $k$
is denoted by $F^k$.
We call $F$ {\em regular\/} if
\begin{itemize}
\item
$F$ is closed under subformulas and renaming of bound variables,
\item
with each $A(x_1,\ldots,x_s)$ in $F$, the class $F$ contains a normal form
of $A$,
\item
for any $k\ge1$, $F^k$ has the {\em pattern set\/}
$P^k\subseteq\{\forall,\exists\}^k$ such that
a normal sentence $A$ belongs to $F^k$ iff every sequence of nested quantifiers
of $A$ belongs to~$P^k$. (By the normality of $A$ all
quantifier sequences have the same length.) 
\end{itemize}

\begin{theorem}\label{thm:lvsdf}
Suppose that $F$ is regular and $G$ is definable in $F$. Then
$$
L_F(G)<\tower\of{D_F(G)+\log^*D_F(G)+2}.
$$
\end{theorem}

\noindent
Note that Theorem \ref{thm:lvsdf} generalizes Lemma \ref{thm:lvsda}
because the classes $\Lambda_a$ are regular.

When we write $\bar z$, we will mean an $s$-tuple $(z_1,\ldots,z_s)$.
If $\bar u\in V(G)^s$, we write $G,\bar u\models A(\bar x)$
if $A(\bar x)$ is true on $G$ with each $x_i$ assigned the respective~$u_i$.
Notation $\models A(\bar x)$ will mean that $A(\bar x)$ is true
on all digraphs with $s$ designated vertices.

A formula $A(\bar x)$ in $F$ is called an {\em $F$-description\/} of
$(G,\bar u)$ if
\begin{itemize}
\item
$G,\bar u\models A(\bar x)$, and
\item
for every $B(\bar x)\in F$ such that $G,\bar u\models B(\bar x)$, we have
$\models A(\bar x)\Rightarrow B(\bar x)$, where $X\Rightarrow Y$ is a shorthand for
$(\neg X)\vee Y$.
\end{itemize}

The next proposition will be useful.

\begin{lemma}\label{lm:descr} 
Suppose that $G$ is definable in $F$. Let $A$ be a sentence in $F$. 
Then $A$ defines $G$ iff
$A$ is an $F$-description of $G$.
\end{lemma}
 \begin{proof} Suppose that $A$ defines $G$. Then $G\models A$. Let 
$B\in F$
satisfy $G\models B$. We have to
show that $H\models A\Rightarrow B$ for any $H$.  
If $H\not\models A$, we are done immediately.  If $H\models A$, then
$H\cong G$ and $H\models B$, as required.

For the other direction, suppose that $A$ is an $F$-description of $G$. We
have to show that $H\not\models A$ for any $H\not\cong G$.
Fix a sentence 
$B\in F$ defining $G$. Since $H\not\models B$ and
$\models A\Rightarrow B$, we conclude that $H\not\models A$, as required.\end{proof}

Let $G$ and $H$ be digraphs, $\bar u\in V(G)^s$, and $\bar v\in V(H)^s$.
We write $G,\bar u\equiv H,\bar v\mod F$ if, for any $A(\bar x)$ in $F$,
we have $G,\bar u\models A(\bar x)$ exactly when $H,\bar v\models A(\bar x)$.

\begin{lemma}\label{lm:descrcl}
Suppose that $G,\bar u\equiv H,\bar v\mod F$ and let $A(\bar x)\in F$.
Then $A$ is an $F$-description of $(G,\bar u)$ iff it is an $F$-description
of $(H,\bar v)$.
\end{lemma}

\begin{proof}
As $A$ is in $F$, we have $G,\bar u\models A(\bar x)$ iff
$H,\bar v\models A(\bar x)$.
Let $B(\bar x)\in F$. Again $G,\bar u\models B(\bar x)$ iff
$H,\bar v\models B(\bar x)$. It follows that $G,\bar u\models B(\bar x)$
implies $\models A(\bar x)\Rightarrow B(\bar x)$ iff $H,\bar v\models B(\bar x)$
implies $\models A(\bar x)\Rightarrow B(\bar x)$.
\end{proof}

Furthermore, we define
$$
(G,\bar u)\bmod F=\setdef{(H,\bar v)}{G,\bar u\equiv H,\bar v\mod F}.
$$
 Let $0\le s\le k$. The class of formulas in $F$ with $s$ free
variables and quantifier rank $k-s$ is denoted by $F^{k,s}$.
In particular, $F^{k,0}=F^k$. We define
$$
E(F^{k,s})=\setdef{(G,\bar u)\bmod F^{k,s}}{G\mbox{\ is\ a\ digraph},\,
\bar u\in V(G)^s}.
$$
We will also use the following notation.
Given $P^k\subseteq\{\forall,\exists\}^k$
and $\gs\in\{\forall,\exists\}^s$, let
$$
P^{k,s}_\gs=\setdef{\rho\in\{\forall,\exists\}^{k-s}}{\gs\rho\in P^k}.
$$
Furthermore, given a regular $F$ with pattern set $P^k$,
let $F^{k,s}_\gs$ consist of the normal formulas in $F^{k,s}$
whose sequences of nested quantifiers are in $P^{k,s}_\gs$. We say that a formula $A(\bar x)\in F^{k,s}_\gs$ {\em describes\/}
a class $\ga\in E(F^{k,s}_\gs)$ if $A(\bar x)$ is an
$F^{k,s}_\gs$-description of some $(G,\bar u)\in\ga$.
By Lemma \ref{lm:descrcl}, this definition does not depend on the particular
choice of a representative $(G,\bar u)$ of $\ga$, and the word 
{\em some\/} in the
definition can be replaced with {\em every}.

\begin{proofof}{Theorem \ref{thm:lvsdf}} For each $\ga\in E(F^{k,s}_\gs)$ we
will construct a formula $A_\ga(\bar x)\in F^{k,s}_\gs$
describing $\ga$. We will use induction on $k-s$. 
Afterwards we will estimate the length of the obtained $A_\ga$ and show how
this implies Theorem \ref{thm:lvsdf}.

We start with $s=k$. Let $\gs\in\{\forall,\exists\}^k$. Assume that $\gs\in
P^k$, for otherwise $F^{k,k}_\gs$ is empty and there is nothing to do.
For any such $\gs$, $F^{k,k}_\gs=F^{k,k}$ is exactly the class of all
quantifier-free formulas in $F$ over the set of variables
$\{x_1,\ldots,x_k\}$. Clearly,  $(G,\bar u)\equiv (H,\bar v) \mod
{F^{k,k}_\gs}$ iff the componentwise 
correspondence between $\bar u$ and $\bar v$ gives a partial isomorphism. So, any given class $\alpha\in
E(F^{k,s}_\gs)$ can be described as follows. Pick any representative
$(G,u_1,\dots,u_k)$ of $\alpha$ and let
$A_\ga(x_1,\ldots,x_k)$ be the conjunction of all atomic formulas
$x_i\mapsto x_j$ for $(u_i,u_j)$ in $G$,
all negations $\neg(x_i\mapsto x_j)$ for $(u_i,u_j)$ not in $G$,
all $x_i=x_j$ for identical $u_i,u_j$, and
all $\neg(x_i=x_j)$ for distinct $u_i,u_j$. Clearly, $H,\bar v\models
A_\ga(\bar x)$ iff $(H,\bar v)\in \alpha$. It follows that 
$A_\ga$ indeed
describes~$\alpha$. Note that $L(A_\ga)\le 18k^2$.

Assume now that $0\le s<k$ and that for any
$\tau\in\{\forall,\exists\}^{s+1}$ with $F^{k,s+1}_\tau\not=\emptyset$ 
and $\gb\in E(F^{k,s+1}_\tau)$ we have
a formula $A_\gb(\bar x,x_{s+1})\in F^{k,s+1}_\tau$ describing~$\gb$.
Given a digraph $G$, an $s$-tuple of vertices $\bar u\in V(G)^s$,
and a non-empty class of formulas $F$, we set
$$
S(G,\bar u;F)=\setdef{(G,\bar u,u)\bmod F}{u\in V(G)}.
$$
We also set $S(G,\bar u;\emptyset)=\emptyset$.
We will write $A\doteq A'$ if formulas $A$ and $A'$ are literally identical.
Let $\gs\in\{\forall,\exists\}^s$ and $\ga\in E(F^{k,s}_\gs)$.
To construct $A_\ga(\bar x)$, we fix $(G,\bar u)$ being an arbitrary
representative of $\ga$ and put%
\footnote{
Here $A_\ga$ has the same form as the {\em Hintikka formula\/}
in~\cite[page 18]{ebbinghaus+flum:fmt}. Curiously, in a similar context in
\cite[Lemma 3.4]{\PSV} we use another generic defining formula
borrowed from \cite[Theorem 2.3.2]{spencer:slrg}, which is not usable now
because $F$ may be not closed with respect to negation.}
$$
A_\alpha(\bar x)\doteq
\bigwedge_{\beta\in S(G,\bar u;F^{k,s+1}_{\gs\exists})}
\exists x_{s+1}\, A_\beta(\bar x,x_{s+1})
\ \wedge\ \forall x_{s+1}
\bigvee_{\beta\in S(G,\bar u;F^{k,s+1}_{\gs\forall})} A_\beta(\bar x,x_{s+1}).
$$

\begin{claim}
$A_\ga(\bar x)\in F^{k,s}_\gs$.
\end{claim}

\begin{subproof}
This follows from 
the assumption that $A_\gb(\bar x,x_{s+1})\in F^{k,s+1}_{\gs{*}}$
for $\gb\in S(G,\bar u;F^{k,s+1}_{\gs{*}})$.
\end{subproof}

\begin{claim}
$G,\bar u\models A_\ga(\bar x)$.
\end{claim}

\begin{subproof}
 Let us show first that all conjunctions over $\beta\in S(G,\bar
 u;F^{k,s+1}_{\gs\exists})$ are satisfied. Each such $\beta$ is of the form
$(G,\bar u,u_\gb)\bmod F^{k,s+1}_{\gs\exists}$ for some $u_\gb\in V(G)$.
By assumption, $G,\bar u,u_\gb\models A_\gb(\bar x,x_{s+1})$ and hence
$G,\bar u\models\exists x_{s+1}A_\gb(\bar x,x_{s+1})$.

It remains to show that the universal member of the conjunction is also
satisfied. Consider an arbitrary $u\in V(G)$. Let
$\gb_u=(G,\bar u,u)\bmod F^{k,s+1}_{\gs\forall}$. By assumption,
$G,\bar u,u\models A_{\gb_u}(\bar x,x_{s+1})$ and hence the disjunction is
always true.
\end{subproof}

\begin{claim}
We have $\models A_\ga(\bar x)\Rightarrow B(\bar x)$ for any $B(\bar x)\in
F^{k,s}_\gs$ such that 
 \beq{B}
 G,\bar u\models B(\bar x).
 \eeq
\end{claim}

\begin{subproof}
Given a class of formulas $F$, let $F|_{\exists}$ (resp.\ $F|_{\forall}$)
denote the class of those formulas in $F$ having form $\exists x (\ldots)$
(resp.\ $\forall x (\ldots)$). First, we settle two special cases of the
claim.

\case 1{$B\in F^{k,s}_\gs|_{\exists}$}
Let $B\doteq\exists x_{s+1} C(\bar x,x_{s+1})$. Note that
$C(\bar x,x_{s+1})\in F^{k,s+1}_{\gs\exists}$.
Assume that $H,\bar v\models A_\ga(\bar x)$. We have to verify that $H,\bar
v\models B(\bar x)$. By~\req{B} we can choose a
vertex $u\in V(G)$  such
that $G,\bar u,u\models C(\bar x,x_{s+1})$. Let
 $$\gb=(G,\bar u,u)\bmod F^{k,s+1}_{\gs\exists}.$$
 We have $H,\bar v\models\exists x_{s+1}A_\gb(\bar x,x_{s+1})$
and hence $H,\bar v,v\models A_\gb(\bar x,x_{s+1})$ for some $v\in V(H)$.
Since we have assumed that $A_\gb(\bar x,x_{s+1})$ is an
$F^{k,s+1}_{\gs\exists}$-description of $\beta$, we have 
$H,\bar v,v\models C(\bar x,x_{s+1})$ and hence
$H,\bar v\models B(\bar x)$ as needed.

\case 2{$B\in F^{k,s}_\gs|_{\forall}$}
 Let $B\doteq\forall x_{s+1} C(\bar x,x_{s+1})$. Note that
$C(\bar x,x_{s+1})\in F^{k,s+1}_{\gs\forall}$.
Assume that $H,\bar v\models A_\ga(\bar x)$. It follows that
for every $v\in V(H)$ there is a $\gb_v\in S(G,\bar u;F^{k,s+1}_{\gs\forall})$
such that $H,\bar v,v\models A_{\gb_v}(\bar x,x_{s+1})$.
By~\req{B} we have $G,\bar u,u\models C(\bar x,x_{s+1})$ for all $u\in V(G)$.
Let $u_v$ be such that 
 $$\gb_v=(G,\bar u,u_v)\bmod F^{k,s+1}_{\gs\forall}.$$
 We have $G,\bar u,u_v\models C(\bar x,x_{s+1})$,
and, by our assumption that $A_{\gb_v}$ describes $\gb_v$, we have
$H,\bar v,v\models C(\bar x,x_{s+1})$. Since $v$ is arbitrary, we conclude
that $H,\bar v\models B(\bar x)$, finishing the proof of Case~2.\medskip

Finally, take an arbitrary $B(\bar x)\in F^{k,s}_\gs$.
Since $B$ is normal (and $s<k$), it is equivalent to a DNF formula
$\vee_i(\wedge_j B_{i,j})$ with all $B_{i,j}$ belonging to
$F^{k,s}_\gs|_{\exists}\cup F^{k,s}_\gs|_\forall$. This can be routinely shown
by induction on $L(B)$. For example, if $B=B_1\wedge B_2$,
where, by induction, $B_h$ is equivalent to 
$\vee_{i_h}(\wedge_{j_h} B_{i_h,j_h})$, $h=1,2$, then we can take
$\vee_{i_1,i_2} ((\wedge_{j_1}
B_{i_1,j_1})\wedge(\wedge_{j_2} B_{i_2,j_2}))$ for $B$. 

Since $G,\bar u\models B(\bar x)$, we have $G,\bar u\models B_{i_0,j}(\bar
x)$ for some $i_0$ and all $j$. From Cases~1--2 it follows that $H,\bar
v\models B_{i_0,j}(\bar x)$ for all $j$ whenever $H,\bar v\models
A_\ga(\bar x)$. This means that $H,\bar v\models B(\bar x)$ whenever $H,\bar
v\models A_\ga(\bar x)$, as required.
\end{subproof}

Let us now estimate the length of the constructed formulas. The estimates are
similar to those in~\cite[Theorem~10.1]{\PSV}. Our bound will
depend on $k$ and $s$ only, so we define
 $$
 l(k,s)=\max_{\tau\in\{\forall,\exists\}^{s}}\max
\setdef{L(A_\gb)}{\gb\in E(F^{k,s}_\tau)}.
 $$
 Let $f(k,s)= |\EHRV{k,s}|$, where $\EHRV{k,s}=E(FO^{k,s})$ with $FO$
being the class of all first order formulas. (The elements of $\EHRV{k,s}$ are
called {\em digraph 
Ehrenfeucht values\/}, see  \cite{spencer:slrg}.) The function
$f(k,s)$ is an upper bound on $|E(F^{k,s}_\tau)|$ for any
$\tau\in\{\forall,\exists\}^{s}$.  The number of Ehrenfeucht values for
(unoriented) graphs was estimated in \cite[Theorem
2.2.1]{spencer:slrg}. The obvious modifications of the proofs from
\cite{spencer:slrg}  give the following bounds for
digraphs: 
 \begin{eqnarray*}
 f(k,k)&\le & 4^{k^2},\\
 f(k,s)&\le & 2^{f(k,s+1)}.
 \end{eqnarray*} 

We already know that $l(k,k)\le18k^2$. 
Also, the analysis of our construction shows
that for $0\le s< k$ we have
 \beq{rec_l}
 l(k,s)\le2f(k,s+1)(l(k,s+1)+9).
 \eeq%
 \comment{Here I am not sure where $+9$ comes from (I get smaller). 
But this constant is not
too important, so I did not bother.}

Let $k\ge 2$. 
Set $g(x)=2\cdot 2^x(x+9)$. A simple inductive argument shows that 
 $$
 f(k,s)\le 2^{g^{(k-s)}(18k^2)}\quad \mbox{and}\quad l(k,s)\le
g^{(k-s)}(18k^2). 
 $$%
 \comment{where $g^{(i)}(x)=g(g(\dots g(x)\dots)$ is obtained by iteratively
applying $g$ $i$ times.}
 Define the two-parameter function 
$\tower(i,x)$ inductively on $i$ 
by $\tower(0,x)=x$ and $\tower(i+1,x)=2^{\tower(i,x)}$ for $i\ge 1$.
This is a
generalization of the old function:
$\tower(i,1)=\tower(i)$. One can prove by induction on $i$ that for any
$x\ge 5$ and $i\ge 1$ we have 
 \beq{g}
 g^{(i)}(x)< \tower(i+1,x)/2.
 \eeq 
 Indeed, it is easy to check the validity of \req{g} for $i=1$, while for
$i\ge 2$ we have
 \beq{indg}
 g^{(i)}(x)< g(\tower(i,x)/2) < 2^{\tower(i,x)-1} = \tower(i+1,x)/2.
 \eeq%
 \comment{Here we used the inequality $g(y/2)\le 2^{y-1}$ valid for all $y\ge
 12$. (The restriction $x\ge 5$ from from $\tower(1,x)/2\ge 12$.}

If $k\ge 12$, then $18k^2< 2^k$ and by~\req{g} we have
 \beq{l}
 l(k,0)\le g^{(k)}(18k^2) < \tower(k+1,18k^2)/2 < \tower(k+\log^*k+2)
 \eeq
 Also, $18\cdot 11^2 < \tower(4)/2$ and, similarly to~\req{indg}, we have 
$g^{(k)}(18k^2)< \tower(k+4)/2$ for $k\le 11$. Thus~\req{l} holds for
$k\in[3,11]$ too. For $k=2$ one can still  prove~\req{l} using~\req{rec_l} and
the sharper initial estimates $f(2,2)= 10$ and $l(2,2)\le 24$.%
 \comment{Here are the calculations:
 \begin{eqnarray*}
 l(2,1)&\le& 2\cdot 10\cdot 33=660,\\
 l(2,0)&\le& 2\cdot 2^{10} \cdot 669.
 \end{eqnarray*}
 The estimate of $f(2,2)$: the given two vertices can span at most $4$ arcs
(we allow loops) and the number of configurations with $0,1,2,3,4$
arcs is respectively $1,2,4,2,1$.}

To finish the proof of Theorem \ref{thm:lvsdf}, let $k=D_F(G)\ge D(G)\ge 2$ 
and $\ga=G\bmod{F^{k,0}}$. Since $G$ is definable in $F^{k,0}$,
the sentence $A_\ga$ defines $G$ by Lemma \ref{lm:descr}.
By~\req{l}, 
 $$L_F(G)\le L(A_\ga)\le l(k,0)< \tower(k+\log^*k + 2),$$ 
 completing the proof.
\end{proofof}

\small

%\bibliographystyle{amsplain}
%\bibliography{oleg,general,graph,misc,random,ramsey}
%\end{document}

\providecommand{\bysame}{\leavevmode\hbox to3em{\hrulefill}\thinspace}
\providecommand{\MR}{\relax\ifhmode\unskip\space\fi MR }
\providecommand{\MRhref}[2]{%
  \href{http://www.ams.org/mathscinet-getitem?mr=#1}{#2}
}
\providecommand{\href}[2]{#2}

\end{document}